\documentclass[11pt,oneside]{article}
\usepackage{amssymb,amsmath,amsthm}
\usepackage{url}
\usepackage[english]{babel}
\usepackage{pstricks,pst-node}
\renewcommand{\le}{\leqslant}
\renewcommand{\ge}{\geqslant}
\newcommand{\bad}{\mathbf{Bad}}
\newcommand{\pbad}{\mathbf{PBad}}

\renewcommand{\AA}{\mathbb{A}}
\newcommand{\RR}{\mathbb{R}}
\newcommand{\CC}{\mathbb{C}}
\newcommand{\ZZ}{\mathbb{Z}}
\newcommand{\QQ}{\mathbb{Q}}
\newcommand{\PP}{\mathbf{P}}

\newcommand{\NN}{\mathbb{N}}

\newcommand{\UU}{\mathbf{U}}
\newcommand{\TT}{\mathrm{T}}

\newcommand{\WWW}{\mathcal{W}}

\newcommand{\vw}{\mathbf{w}}
\newcommand{\vx}{\mathbf{x}}
\newcommand{\vu}{\mathbf{u}}
\newcommand{\vv}{\mathbf{v}}

\newcommand{\vq}{\mathbf{q}}

\newcommand{\vd}{\mathbf{d}}

\newtheorem{lemma}{Lemma}
\newtheorem{theorem}{Theorem}
\newtheorem*{ntheorem}{Theorem}
\newtheorem{proposition}{Proposition}
\newtheorem{problem}{Problem}

\newtheorem*{conjecture}{Conjecture}

\newtheorem*{corollary}{Corollary}
\newtheorem{corol}{Corollary}
\newtheorem{fact}{Fact}

\newcommand{\Sl}{\mathrm{SL}}
\newcommand{\Gl}{\mathrm{GL}}
\newcommand{\mad}{\mathbf{Mad}}
\newcommand{\lmad}{\mathbf{LMad}}
\newcommand{\tr}{\mathrm{tr}}
\newcommand{\Id}{\mathrm{Id}}

\begin{document}

\title{On continued fraction expansion of potential counterexamples to $p$-adic Littlewood
conjecture}

\author{
 Dzmitry Badziahin
\footnote{Research supported by EPSRC  Grant EP/L005204/1} }

\maketitle

\begin{abstract}
The $p$-adic Littlewood conjecture (PLC) states that
$\liminf_{q\to\infty} q\cdot |q|_p \cdot ||qx|| = 0$ for every prime
$p$ and every real $x$. Let $w_{CF}(x)$ be an infinite word composed
of the continued fraction expansion of $x$ and let $\TT$ be the
standard left shift map. Assuming that $x$ is a counterexample to
PLC we show that limit elements of the sequence $\{\TT^n
w_{CF}(x)\}_{n\in\NN}$ are quite natural objects to investigate in
attempt to attack PLC for $x$. We then get several quite restrictive
conditions on such limit elements~$w$. As a consequence we prove
that we must have $\lim_{n\to\infty} P(w,n) - n = \infty$ where
$P(w,n)$ is a word complexity of $w$. We also show that $w$ can not
be among a certain collection of recursively constructed words.
\end{abstract}

\section{Introduction}

In 2004 de~Mathan and Teulie~\cite{Mathan_Teulie_2004} proposed the
following problem which is now called a $p$-adic Littlewood
conjecture (PLC).
\begin{conjecture}[PLC]
Let $p$ be a prime number. Then every real $x$ satisfies
\begin{equation}\label{eq_main}
\liminf_{q\to \infty} q\cdot |q|_p\cdot ||qx|| =0
\end{equation}
where $||\cdot||$ denotes the distance to the nearest integer.
\end{conjecture}

It is widely believed to be easier than the famous Littlewood
conjecture where the expression above is replaced by
$$
\liminf_{q\to \infty} q\cdot ||qx||\cdot ||qy|| =0.
$$
Despite of the essential efforts from the mathematical communities
both conjectures still remain open.

Assume that there is a counterexample $x$ to PLC. Then it must
satisfy
\begin{equation}\label{cond_counter}
\inf_{q\in\NN} \{q \cdot |q|_p\cdot ||qx||\}\ge \epsilon
\end{equation}
for some $\epsilon>0$. We denote the set of $x\in\RR$ which satisfy
the condition~\eqref{cond_counter} by $\mad_\epsilon$. So PLC is
equivalent to saying that for all $\epsilon>0$ the set
$\mad_\epsilon$ is empty.

It is already known that $\mad_\epsilon$ is very ``small''. The
condition~\eqref{cond_counter} straightforwardly implies that $x$ is
badly approximable hence it belongs to the set of zero Lebesgue
measure. Moreover $\mad_\epsilon$ is included in the subset
$\bad_\epsilon$ of the set $\bad$ of badly approximable numbers
which is defined as follows
$$
\bad_\epsilon:=\{x\in \RR\;:\; \inf_{q\in\NN}\{q\, ||qx||\}\ge
\epsilon\}.
$$
In~\cite{Mathan_Teulie_2004} it was shown that quadratic
irrational~$x$, the classical examples of badly approximable
numbers, satisfy PLC. Later Bugeaud, Drmota and
de~Mathan~\cite{Bugeaud_Drmota_Mathan_2007} generalized that result
to numbers which continued fraction expansion contain arbitrarily
long periodic blocks. In 2007 Einsiedler and
Kleinbock~\cite{Einsiedler_Kleinbock_2007} proved that
$\mad_\epsilon$ is of zero box dimension for every $\epsilon>0$.

Numbers from $\bad$ can easily be described in terms of their
continued fraction expansion. We'll state this classical fact in
terms of infinite words (the details can be found
in~\cite{Cassels_1955}, for example).

\begin{fact}\label{fact1}
Let $w_{CF}(x)$ be an infinite word composed of the partial
quotients of the continued fraction expansion of
$x=[0;a_1,a_2,\ldots]\in\RR$. If $x\in\mad_\epsilon$ then
$w_{CF}(x)\in \AA_N^\NN$ where $N = \epsilon^{-1}+1$ and $\AA_N:=
\{1,2,\ldots,N\}$. In the other words the word $w_{CF}(x)$ belongs
to the finite alphabet $\AA_N$.
\end{fact}

So if one wants to attack the PLC then it is natural to consider
numbers $x$ with $w_{CF}(x)\in\AA_N^\NN$ and to try to impose as may
conditions on $w_{CF}(x)$ for potential counterexamples to PLC as
possible. The ideal situation of course would be if that conditions
become self-contradictory which would immediately imply that $\mad =
\emptyset$.

Some restrictions on $w_{CF}(x)$ for potential counterexamples $x$
to PLC were recently imposed in~\cite{BBEK_2014}. In order to state
them we introduce a couple of new definitions.
We call a word $w\in \NN^\NN$ recurrent if every finite block
occurring in $w$ occurs infinitely often. Then $w$ is eventually
recurrent if $\TT^m w$ is recurrent for some positive integer $m$
where $\TT$ is a standard left shift map on $\NN^\NN$.

\begin{ntheorem}[BBEK1]
If $w_{CF}(x)$ is eventually recurrent then $x$ satisfies PLC.
\end{ntheorem}

Almost all well known classical infinite words such as Sturmian
words or the Thue-Morse word, are recurrent so Theorem~BBEK1 states
that none of them can be a continued fraction expansion of the
counterexample to PLC. Another implication of this theorem stated
in~\cite{BBEK_2014} is that for $x\in\mad$ the word complexity
$P(w_{CF(x)},n)$ as a function of $n$ does not grow too slow. By the
complexity $P(w,n)$ of the word $w$ we mean the number of distinct
blocks of length $n$ which occur in $w$.

\begin{corollary}[BBEK2]
If $x\in \mad$ then $P(w_{CF}(x),n) - n\to\infty$ as $n\to\infty$.
\end{corollary}

On the other hand the next result from~\cite{BBEK_2014} states that
the complexity of $x\in\mad$ can not grow too fast as well.

\begin{ntheorem}[BBEK3]
If $x\in\mad$ then
$$
\lim_{n\to\infty} \frac{\log P(w_{CF}(x), n)}{n} = 0.
$$
\end{ntheorem}
\noindent In other words Theorem BBEK3 says that the complexity
$w_{CF}(x)$ of a counterexample to PLC grows slower than any
exponential function. Surely there is a huge gap between
subexponential growth from Theorem BBEK3 and linear growth from
Corollary BBEK2. So there are still plenty of words $w_{CF}(x)$
uncovered by both of these statements.

In this paper we put more restrictions on the continued fraction
expansion of the counterexamples to PLC. The next section explains
that it is quite reasonable to look at the limit points of $\{\TT^n
w_{CF}(x)\}_{n\in\NN}$ in order to attack PLC. For brevity the set
of limit points of $\{\TT^n w_{CF}(x)\}_{n\in\NN}$ will be called by
$\WWW_{CF}(x)$. In particular we provide several ``easy-to-state''
conditions on $w\in \WWW_{CF}(x)$ which imply PLC for $x$. One of
them (Theorem~\ref{th_sturmian}), for example, substantially
improves Corollary~BBEK2 by putting the complexity condition from
this corollary on every limit word $w\in\WWW_{CF}(x)$.
Section~\ref{sec_res} introduces a couple of quite general and
rather technical core results on words $w\in\WWW_{CF}(x)$ which are
then proved in the next Section~\ref{sec_mainproof}. The rest of the
paper is devoted to the proofs of all their applications.


\section{A version of PLC for limit words $w\in\WWW_{CF}(x)$. Sets $\lmad$ and $\lmad_\epsilon$}

One may try to attack PLC by taking $q$ to be a linear combination
of the denominators $q_n$ and $q_{n+1}$ of two consecutive
convergents to $x$. The result of this attempt is presented in the
following

\begin{proposition}\label{lem1}
Let $q_n< q_{n+1}$ be the denominators of two consecutive
convergents to $x\in\bad_\epsilon$. Then for every $a,b\in\ZZ$ and
$r = r_{a,b,n} =|aq_n + bq_{n+1}|$ we have
\begin{equation}\label{lem1_eq}
r\cdot |r|_p\cdot ||rx||\le 4\max\{a^2,b^2\}(N+1)\cdot |r|_p
\end{equation}
where $N = [\epsilon^{-1}]+1$.
\end{proposition}

\proof We just use two standard facts: $(N+1) q_n> q_{n+1}$ (by
Fact~\ref{fact1}) and $q_n||q_nx||< 1$, $q_{n+1}||q_{n+1}x||< 1$.
Then
$$
r\cdot |r|_p\cdot||rx||\le (2\max\{|a|,|b|\})^2\cdot q_{n+1} \cdot
\frac{1}{q_n}\cdot |r|_p \le 4\max\{a^2,b^2\}(N+1)|r|_p.\\[-4ex]
$$
\endproof
As the consequence of the proposition, if $x$ is in $\mad_\epsilon$
then for every pair $\vq_n = (q_n,q_{n+1})$ and $a,b\in \ZZ$ one
must have
\begin{equation}\label{eq_pair}
r_{a,b,n}= 0 \mbox{ or }|r_{a,b,n}|_p\ge\frac{\epsilon}{4(N+1)}
\cdot \min\{a^{-2},b^{-2}\}.
\end{equation}
This fact imposes an additional condition on $x\in \mad_\epsilon$ or
more exactly on every pair of denominators of consecutive
convergents to $x$. Every such a pair must satisfy~\eqref{eq_pair}.
By combining this with Fact~1  we get

{\it if $x\in \mad_\epsilon$ then for every $n\in \NN$, $\TT^n
w_{cf}(x)\in
\AA_N^\NN$ and $(q_n,q_{n+1})$ satisfies~\eqref{eq_pair}.}\\[-1ex]

For convenience we introduce the following notation
\begin{equation}\label{eq_aw}
\mathrm{for }\; a\in\NN,\;
A_a:=\left(\begin{array}{cc}0&1\\1&a\end{array}\right);\quad
\mathrm{for }\;w = a_1\ldots a_n,\; A_w = A_{a_1}\cdot \ldots\cdot
A_{a_n}.
\end{equation}
Classical results from the theory of continued fractions relate the
pairs $(q_n,q_{n+1})$ for different values of $n$ in the following
way:
$$
\left(\begin{array}{cc}0&1\\1&a_{n+1}\end{array}\right)
\left({q_{n-1}}\atop{q_n}\right) = \left(q_n\atop{q_{n+1}}\right).
$$
This observation gives rise to the following extension of the left
shift map $\TT$ on $\AA_N^\NN$ to the pairs $(w,\vq)\in
\AA_N^\NN\times \PP_{\QQ_p}^1$:
$$
\TT(w,\vq):= (\TT w, A_{a_1} \vq).
$$
Condition~\eqref{eq_pair} is reflected in the following definition.
A point $\vq=(q_1,q_2)\in \PP_{\QQ_p}^1$ is said to be $p$-adically
badly approximable if there exists $\epsilon>0$ such that $\forall
(a,b)\in \ZZ\backslash \{(0,0)\}$ one has
$$
|aq_1 +bq_2|_p\cdot \min\{|q_1^{-1}|_p,|q_2^{-1}|_p\}\ge
\min\{|a|^{-2},|b|^{-2}\} \cdot \epsilon.
$$
Sometimes instead of projective coordinates we will use affine ones.
In that case we say that $\omega\in\QQ_p$ is $p$-adically badly
approximable if $(\omega,1)$ is. We call the set of all $p$-adically
badly approximable points by $\pbad$. Then by analogy with the
definition of $\bad_\epsilon$ we define the set $\pbad_\epsilon$ as
the subset of $\pbad$ containing those points $\vw$ for which
$\epsilon$ in the definition is fixed.

The next crucial proposition links the counterexamples $x$ to PLC
with the limit points of the sequence $\{T^n(w_{CF}(x),({0\atop
1}))\}_{n\in\NN}$.

\begin{proposition}\label{corl_t}
If $x$ is a counterexample to PLC then every limit point $(w,\vx_p)$
of $\{T^n(w_{CF}(x),({0\atop 1}))\}_{n\in\NN}$ satisfies the
following property: there exists $\epsilon>0$ such that for every
$n\in\NN$ one has $T^n(w,\vx_p)\in \AA_{N}^\NN\times \pbad_\epsilon$
where $N = [\epsilon^{-1}]+1$ as above.
\end{proposition}

The proof of the proposition is postponed till the end of this
section. We define the set of pairs $(w,\vx_p)$ which satisfy the
conditions of this proposition by $\lmad$ (letter L is for limit
point). By $\lmad_\epsilon$ we define the subset of $\lmad$
containing those pairs for which the parameter $\epsilon$ is fixed.
In other words the proposition says that if at least one limit point
of $\{T^n(w_{CF}(x),({0\atop 1}))\}_{n\in\NN}$ is not in $\lmad$
then $x$ is not in $\mad$. By Fact~\ref{fact1} and the compactness
of the set $\AA^\NN_N$ for every $x\in\bad_\epsilon$ elements
$\TT^n(w_{CF}(x))$ belong to the compact set. It is also well known
that $\PP_{\QQ_p}^1$ is compact. Therefore the sequence
$\{T^n(w_{CF}(x),({0\atop 1}))\}_{n\in\NN}$ must have at least one
limit point which can be tested for inclusion in $\lmad$.

The upshot of this discussion is that to attack PLC it is natural to
investigate the set $\lmad$. The straightforward problem about it is
\begin{problem}
Is $\lmad$ empty?
\end{problem}
Surely the positive answer to this problem would immediately imply
that $\mad$ is empty as well. Unfortunately this is not the case as
shown by the following result. Let $w=a_1a_2\ldots$ be an infinite
word. Then by $w_l$ we denote its finite prefix of length $l$:
$w_l:= a_1\ldots a_l$.
\begin{theorem}\label{th_lmad}
Let $w\in\AA_N^{\NN}$ be a periodic infinite word and $l\in\NN$ be
the length of its minimal period. Then $(w,\vx_p)$ is in $\lmad$ if
and only if $\vx_p$ coincides with one of the eigenvectors
of~$A_{w_l}^T$.
\end{theorem}

Since infinitely (countably) many different matrices $A_w^T$ have
eigenvectors in $\PP_{\QQ_p}^1$, this theorem shows that the set
$\lmad$ is infinite. And therefore the answer to Problem A is
``no''. Luckily or not, Theorem~\ref{th_lmad} does not provide us
with any counterexample to PLC. Indeed, Bugeaud, Drmota and de
Mathan~\cite{Bugeaud_Drmota_Mathan_2007} showed that for $x\in\mad$
a periodic word~$w$ can not be in $\WWW_{CF}(x)$.

We will prove Theorem~\ref{th_lmad} in Section~\ref{sec_5}.

\medskip\noindent{\bf Remark.} In~\cite{Einsiedler_Kleinbock_2007}
the following generalization of PLC was posed: {\it every pair
$(x,y)\in {\RR_{>0}}\times {{\QQ}_p}$ satisfies}
$$
\inf_{a\in\NN,b\in \NN\cup\{0\}} \max\{|a|,|b|\}\cdot |ax-b|\cdot
|ay-b|_p =0. \eqno{(EK)}
$$
This conjecture has close connection with Problem~A. In particular a
similar method to that used in the proof of Theorem~\ref{th_lmad}
shows that if $x\in\RR_{>0}$ and $y\in\QQ_p$ are irrational roots of
the same quadratic polynomial then Condition~(EK) fails which in
turn means that the proposed generalization of PLC is false. We
leave the details of this fact to the interested reader as an
exercise.

The next natural question about $\lmad$ is: are there more elements
in it? Or more formally we state it as follows
\begin{problem}
Is it true that for every $(w,\vx_p)\in \lmad$ the word $w$ is
periodic?
\end{problem}
The positive answer to this problem will again imply PLC. We will
reformulate it in a slightly different way which, in authors
opinion, is more natural. It follows straightforwardly from the
definition that the set $\lmad$ is invariant under~$\TT$ and so is
the set $\lmad_\epsilon$ for every $\epsilon>0$. Moreover it can be
verified that $\lmad_\epsilon$ is also closed and therefore it is
compact. It is well known that every compact invariant set contains
a minimal invariant subset i.e. a subset which does not contain any
other non-empty invariant closed subsets. Therefore it is sufficient
to consider the minimal compact invariant subsets of $\lmad$.

Note that for each $\vx\in\lmad$ described in Theorem~\ref{th_lmad}
we have $\TT^l\vx = \vx$ which in turn implies that the set $\{\TT^n
\vx\}_{n\in\NN}$ is finite. So all minimal invariant sets generated
by these elements $\vx$ are also finite. On the other hand every
finite minimal invariant set must be generated by an element
$(w,\vx_p)$ with periodic $w$ (by finiteness we must have $\TT^l w =
w$ for some $l\in\NN$). Therefore Problem~B can be reformulated as
follows:

\medskip \noindent{\bf Problem B'}\;{\it Are there infinite closed subsets of $\lmad$ which are minimal
invariant under~$\TT$?}

An advantage of investigating minimal invariant subsets is that we
know quite well how should the first coordinate of every element
$(w,\vx_p)$ of a minimal invariant set look like
(\cite{Lothaire_2001}[Theorem 1.5.9]). In that case $w$ is uniformly
recurrent or in other words it is recurrent and for each finite
factor $u$ of $w$ the distance between any two consecutive
appearances of $u$ in $w$ is bounded above by some constant $d =
d(u)$.

The author believes that Problem~B (respectively~B') has a positive
answer. In this paper we will impose rather restrictive conditions
on the elements of $\lmad$. The most general but technical of them
are presented in Theorems~\ref{th_main} and~\ref{th_da} which are
introduced in Section~\ref{sec_res}. The author does not know if
they rule out every possible infinite word $w$ (this would actually
imply PLC) however it seems that this is not the case, there are
still plenty of uncovered infinite words and therefore there are
many numbers $x$ for which the PLC still remains open. Here we
provide some ``easy-to-state'' applications of that results.

First of all, if $(w,\vx_p)\in\lmad$ then the complexity of the word
$w$ can not grow too slow.
\begin{theorem}\label{th_sturmian}
Let non-periodic $w\in\AA_N^\NN$ be such that $\forall n\in\NN$,
$P(w,n)\le n+C$ for some positive absolute constant $C$. Then $\vx =
(w,\vx_p)\not\in\lmad$ for every $\vx_p\in\PP_{\QQ_p}^1$.
\end{theorem}
The proof of this theorem is described in
Section~\ref{sec_lowcomplex}. Since the complexity function is
strictly increasing then the straightforward corollary of
Theorem~\ref{th_sturmian} is
\begin{corollary}
If $x$ is a counterexample to PLC then for any $w\in \WWW_{CF}(x)$,
$$
\lim_{n\to \infty} P(w,n) - n \to \infty.
$$
In particular, the set $\WWW_{CF}(x)$ must not contain Sturmian
words.
\end{corollary}

Next, we can show that for a big collection of words $w$ which can
be recurrently constructed by concatenations the pair $(w,\vx_p)$ is
never in $\lmad$. We call $\WWW(\sigma_1,\sigma_2,\ldots,\sigma_n)$
a concatenation map if it is some composition of concatenations of
words $\sigma_1,\ldots,\sigma_n$. Recall that by $A_w$ we denote the
matrices defined by~\eqref{eq_aw}. Also for any ring $R$ with
identity by $\Sl_2^\pm (R)$ we denote the set of $2\times 2$
matrices with determinant $\pm1$. Finally by $\ZZ_p$ we denote the
set of $p$-adic integer numbers.

\begin{theorem}\label{th_concat}
Let a sequence of finite words over the alphabet $\AA_N$ be
constructed recursively as follows: $\sigma_{1},\ldots,\sigma_{m}$
are given words of length 1 such that not all of them equal to each
other; for every $n\in\NN$,
$$
\sigma_{n+m}:=
\sigma_{n+m-1}\WWW(\sigma_{n},\sigma_{n+1},\ldots,\sigma_{n+m-1})
$$
where $\WWW$ is a concatenation map. Assume that for every $m$-tuple
of words $\eta_1,\ldots,\eta_m$ the equation
$$
A_{\eta_m} = \widetilde{\WWW}(X, A_{\eta_1},\ldots,A_{\eta_{m-1}})
$$
has at most one solution $X\in \Sl^\pm_2(\ZZ_p)$. Here
$\widetilde{\WWW}$ is made of $\WWW$ by replacing each concatenation
with a product of matrices. Then for every limit word $w$ of the
sequence $\sigma_{n}$ and every $\vx_p\in \PP_{\QQ_p}^1$,
$(w,\vx_p)\not\in \lmad$.
\end{theorem}

In fact the condition on $\sigma_1,\ldots, \sigma_m$ in the theorem
can be weakened. We just need that not all of these words are the
powers of the same finite word. However for the sake of simplicity
we do not put this condition to the theorem. Its proof is provided
in Section~\ref{sec_concat}.

Theorem~\ref{th_concat} covers a big collection of automatic words.
In particular one can check that the Fibonacci word $w_{fib}$
satisfies all of theorem's conditions. Indeed it is the limit point
of the sequence $\{\sigma_n\}_{n\in\NN}$ constructed as follows:
$$
\sigma_1,\sigma_2 \mbox{ are distinct one-digit words;}\quad
\sigma_{n+1} = \sigma_n\sigma_{n-1}.
$$
Therefore $(w_{fib},\vx_p)$ is never in $\lmad$.

Finally we give a nice combinatorial condition on a word $w$ which
guarantees that $(w,\vx_p)$ does never belong to $\lmad$. Before we
do that we define a sequence of bipartite graphs $G_n (S_n, T_n,
E_n)$ related to $w$. Their definition distantly resembles more
classical Rauzy graphs for infinite words. We set both $S_n$ and
$T_n$ to be sets of all different factors of $w$ of length $n$. Then
vertices $s\in S_n$ and $t\in T_n$ are linked with an edge iff the
word $st$ is itself a factor of $w$.

\begin{theorem}\label{th_graph}
Let $w\in \AA_N^\NN$ be a non-periodic uniformly recurrent word. If
a number of connected components in $G_n$ is bounded by an absolute
constant independent of $n$ then $(w,\vx_p)\not\in\lmad$ for every
$\vx_p\in\PP_{\QQ_p}^1$.
\end{theorem}

Note that the number of edges in $G_n(S_n,T_n,E_n)$ coincides with
$P(w,2n)$. Heuristically more edges $G_n$ has, less chance that it
has many connected components is. So it would be natural to test the
conditions of Theorem~\ref{th_graph} for words $w$ of high
complexity.

The proof of this theorem is provided in Section~\ref{sec_graph}.

\noindent{\bf Remark.} None of the results in this section cover the
Thue-Morse word $w_{tm}$. The author believes that the main
Theorem~\ref{th_main} itself can be applied to it in order to show
that $(w_{tm},\vx_p)$ is never in $\lmad$. Anyway it would be
interesting to check if there is any prime $p$ and point $\vx_p\in
\PP_{\QQ_p}^1$ such that $(w_{tm}, \vx_p)\in\lmad$.\\

We finish this section with the proof of Proposition~\ref{corl_t}.
Assume that $x$ is a counterexample to PLC, i.e. $x\in\mad_\delta$
for some $\delta>0$. Choose $\epsilon>0$ such that $\delta =
\epsilon/4(N+1)$. Note that $N_\delta = [\delta^{-1}]+1$ is smaller
than $N=N_\epsilon$ and therefore Fact~1 shows that
$w_{CF}(x)\in\AA_N^\NN$ and consequently every limit point
$w\in\WWW_{CF}(x)$ lies in the same space. Next, notice that $(0,1)$
is a pair $(q_{-1},q_0)$ and therefore
$$
\TT^n\left(w_{CF}(x),\left({0\atop 1}\right)\right) = \left(\TT^n
w_{CF}(x), \left({q_{n-1}\atop q_n}\right)\right).
$$

Consider an arbitrary limit element $(w,\vx_p)$ of
$\{T^n(w_{CF}(x),({0\atop 1}))\}_{n\in\NN}$. Assume that for some
$a,b\in \ZZ$ the value $r_{a,b} = ax_p + bx_p'$ satisfies
$r_{a,b}=0$. It means that there is a subsequence $q_{n_k},
q_{n_k+1}$ for which
$$
r_{a,b,n_k}:= |aq_{n_k} + bq_{n_k+1}|_p
$$
tends to zero as $k\to\infty$. Surely $r_{a,b,n_k}$ can not be zero
for infinitely many of $k$ since the denominators $q_n$ tend to
infinity with $n$. Therefore at some point the expression
$r_{a,b,n_k}$ starts violating~\eqref{eq_pair} and we get a
contradiction. To finish the proof we notice that
Condition~\eqref{eq_pair} remains true for every limit point of the
sequence $(q_n,q_{n+1})_{n\in\NN}$ and therefore it must be true for
$\vx_p$ which immediately implies that $\vx_p\in
\pbad_{\delta/(4(N_\delta+1))}=\pbad_\epsilon$. This finishes the
proof of the proposition.
\endproof

%

\section{Core results}\label{sec_res}

%
%

Before introducing the main result we need to do some preparation.
We equip the space $\PP_{\QQ_p}^1$ with the metrics defined as
follows. Given
$\vw=(\omega_1,\omega_2),\vv=(\upsilon_1,\upsilon_2)\in\PP_{\QQ_p}^1$
let
$$
d(\vw,\vv):=|\omega_1\upsilon_2 - \omega_2\upsilon_1|_p\cdot
\min\{|\omega_1|_p^{-1},|\omega_2|_p^{-1}\}\cdot
\min\{|\upsilon_1|_p^{-1},|\upsilon_2|_p^{-1}\}.
$$
Wherever possible we will take the coordinates $(\omega_1,\omega_2)$
of $\vw$ such that $\max\{|\omega_1|_p,|\omega_2|_p\} = 1$. If this
happens we emphasize it by using slightly different notation:
$\vw\in\tilde{\PP}^{1}_{\QQ_p}$. One can see that for
$\vw,\vv\in\tilde{\PP}^1_{\QQ_p}$ the definition of distance
$d(\vw,\vv)$ as well as the definition of $\vw$ being $p$-adically
badly approximable point becomes simpler.

Consider the set $\Sl^\pm_2(\ZZ_p)$. In accordance to the Jordan
normal form of the matrices in $\Sl^\pm_2(\ZZ_p)$ one can split it
into two subsets:
\begin{itemize}
\item $\Sl^\pm_{2,1}(\ZZ_p)$ which consists of all matrices
having two different eigenvectors in $\PP_{\overline{\QQ}_p}^1$.
\item $\Sl^\pm_{2,2}(\ZZ_p)$ which consists of matrices
which are similar to one of the matrices
$$
\left( {v\;\;0\atop 1\;\;v}\right);\quad v\in\overline{\QQ_p}.
$$
\end{itemize}

\begin{lemma}\label{fact_4}
Let $A\in \Sl^\pm_{2,1}(\ZZ_p)$. Then there exists a positive
integer $\kappa(A)\le p^2$ such that for each eigenvalue $\lambda$
of $A$ one has $|\lambda^{\kappa(A)}-1|_p\le p^{-1}$.
\end{lemma}

\proof Both eigenvalues $\lambda$ and $\lambda'$ of $A$ are the
roots of the quadratic equation
$$
\lambda^2 - \tr(A) \lambda + \det(A)=0.
$$
Since $\det(A) = \pm 1$, $\lambda$ is a unit, so $|\lambda|_p=1$.
Further arguments follow the standard proof of Fermat Little
theorem. We compose a complete list of representatives of residue
classes in $\ZZ_p[\lambda]$ modulo $p$ such that their $p$-adic norm
equals one. Since all that residue classes are contained in the set
$\{ a+b\lambda\;:\; 0\le a,b<p\}$ there are at most $p^2$ elements
on the list. We denote its size by $\kappa(A)$, then the list looks
as follows: $x_1,x_2,\ldots, x_{\kappa(A)}$. If we multiply each
element from the list by $\lambda$ the resulting numbers will again
represent each residue class $\ZZ_p[\lambda]$ modulo $p$ with
$p$-adic norm equals one. Then by looking at products we get
$$
\prod_{i=1}^{\kappa(A)} x_i \equiv \prod_{i=1}^{\kappa(A)} \lambda
x_i \pmod p\quad \Rightarrow\quad |\lambda^{\kappa(A)}-1|_p\le
p^{-1}.
$$
\endproof

Next, we need the notion of the $p$-adic logarithm which is defined
as the following series
$$
\log x = \sum_{n=1}^\infty (-1)^{n-1}\frac{(x-1)^n}{n}.
$$
It is well-defined on disc $|x-1|_p<1$. Therefore by
Lemma~\ref{fact_4} the value $\log \lambda^{\kappa(A)}$ is well
defined for every eigenvalue of $A\in \Sl^\pm_{2,1}(\ZZ_p)$. Notice
also that for $|x|_p\le p^{-1}$ one has $|\log x|_p = |x-1|_p$.
Moreover one can check that for $x,y\in\QQ_p$ such that $|x-1|_p<1$
and $|y-1|_p<1$ the following formula is satisfied
\begin{equation}\label{eq_log}
|\log x-\log y|_p = \left|\log\frac{x}{y}\right|_p =
\left|\frac{x}{y}-1\right|_p = |x-y|_p.
\end{equation}

Now we want to remove (for the reason which will become clear later)
all matrices $A$ from $\Sl^\pm_{2,1}(\ZZ_p)$ such that one of their
eigenvalues $\lambda$ satisfies $\log \lambda^{\kappa(A)}=0$. So we
introduce another set
$$
\widetilde{\Sl}_2(\ZZ_p):=\!\!\left\{A\in
\Sl^\pm_{2,1}(\ZZ_p):\mbox{ for two eigenvectors
$\lambda_1,\lambda_2$ of $A$, }\log \lambda_{1}^{\kappa(A)}\neq 0,
\log\lambda_2^{\kappa(A)}\neq 0\! \right\}.
$$

Given $\vx = (w,\vx_p)\in \AA_N^\NN\times \PP_{\QQ_p}^1$ and
$n\in\NN$ we denote by $\vx_{p,n}$ the element of $\PP_{\QQ_p}^1$
defined by the equality $(T^nw, \vx_{p,n})$. Define
$$
B_k(\vx):=\left\{\left(\alpha\atop \beta\right)\in
\PP_{\QQ_p}^1\;:\; \exists\, n\in\NN,\mbox{ s.t. }
d\left(\left(\alpha \atop \beta\right),\vx_{p,n}\right)\le
p^{-k}\right\}.
$$
Roughly speaking the set $B_k(\vx)$ is comprised by
$p^{-k}$-neighborhoods of the elements $\vx_{p,n}$.
The following theorem lies at the heart of this paper.
\begin{theorem}\label{th_main}
Let $A\in \widetilde{\Sl}_2(\ZZ_p)$. Let $\vw_1 =
\left(\omega_1\atop \omega_2\right)$ and $\vw_2 =
\left(\omega_3\atop \omega_4\right)$ be two eigenvectors of $A^T$
with $\omega_{i}\in \overline{\QQ_p}, i\in\{1,2,3,4\}$ and
$\max\{|\omega_1|_p,|\omega_2|_p\}
=\max\{|\omega_3|_p,|\omega_4|_p\}=1$. Consider $\vx =
(w,\vx_p)\in\AA_N^\NN\times \PP_{\QQ_p}^1$ and define
$$
\epsilon_1:= d(\vx_p,\vw_1),\quad \epsilon_2:=d(\vx_p,\vw_2),
\quad\epsilon_3:= |\log (\lambda_1^{4\kappa(A)})|_p,\quad
\delta:=\min\{\epsilon_1,\epsilon_2,p^{-1}\}
$$
where $\lambda_1$ is one of the eigenvalues of $A$. Assume that
$\delta\neq 0$. Finally let $m\in\NN$ and $k\in\NN$ satisfy the
inequality
\begin{equation}\label{eq_th_1}
p^k\ge \frac{\sqrt{d(\vw_1,\vw_2)\cdot m}}{\epsilon_3\cdot
\delta\cdot \sqrt{2\epsilon_1\epsilon_2\cdot p\cdot \kappa(A)}}.
\end{equation}
If
\begin{equation}\label{eq_th_2}
\{\vx_p, A^T \vx_p , \ldots, (A^T)^m \vx_p\} \subset B_k(\vx)
\end{equation}
then $\vx\not\in \lmad_\epsilon$ for
\begin{equation}\label{eq_epsilon}
\epsilon = \sqrt{\frac{2\epsilon_1\epsilon_2\cdot p\cdot
\kappa(A)}{\epsilon_3^2\delta^2\cdot d(\vw_1,\vw_2)\cdot m}}.
\end{equation}
\end{theorem}

Note that the values $\kappa(A)$, $d(\vw_1,\vw_2)$ and $\epsilon_3$
are solely defined by the matrix $A$ and since $A\in
\widetilde{\Sl}_2(\ZZ_p)$ the last two of them are always strictly
positive. The value $\delta$ measures how ``close'' is $\vx_p$ to
one of the eigenvectors of $A^T$. If the parameters $A$ and $m$ are
fixed then as~$\delta$ tends to zero, the
estimate~\eqref{eq_epsilon} on $\epsilon$ tends to infinity. In
other words smaller the value of $\delta$, weaker the estimate for
$\epsilon$.

Before we move to the next result we show that Theorem~\ref{th_main}
covers all matrices of the form $A_w$ where $w$ is any finite word.

\begin{lemma}\label{fact_6}
For every finite nonempty word $w\in \NN^n$ one has $A_w \in
\widetilde{\Sl}_2(\ZZ_p)$.
\end{lemma}

\proof Firstly we show that $A_w$ has two distinct eigenvalues, so
it is in $Sl_{2,1}^\pm(\ZZ_p)$. The eigenvalues of $A$ are the roots
of the equation $\lambda^2 - \tr(A)\lambda + \det(A) = 0$. Since
$\det(A) = \pm1$ there is only one possibility for $\lambda_1 =
\lambda_2$, namely there should be $\tr(A_w) = 2$ and $\det(A_w) =
1$ but one can easily check that there are no words $w$ with this
property.

For algebraic $\lambda$ the condition $\log \lambda^ {\kappa(A)} =
0$ is only possible when $\lambda$ is a root of unity. The only
quadratic roots of unity are either roots of the equation $x^6=1$ or
$x^4=1$. Therefore they must be roots of one of the following
quadratic equations:
$$
\begin{array}{rrr}
x^2 - 1=0;& x^2 + 1=0;& x^2-2x+1 = 0;\\
x^2+2x+1 = 0;& x^2-x+1 = 0;& x^2+x+1=0.
\end{array}
$$
There are only two words $w_1 = 1$ and $w_2 = 2$ for which $|\tr
(A_w)|\le 2$. However an easy check shows that the eigenvalues of
both $A_1$ and $A_2$ are not roots of unity.\endproof

Now consider different type of matrices $A\in \Sl_2^\pm(\ZZ_p)$
which will be needed later. Namely,
$$
A = D_a :=\left({1\;\; 0\atop a\;\; 1}\right).
$$
It is easily verified that $D_a\not\in\widetilde{\Sl}_2(\ZZ_p)$ so
Theorem~\ref{th_main} is not applicable to it. However very similar
(and even simpler) result is true for such matrices too.
\begin{theorem}\label{th_da}
Let $A = D_a$ for some $a\in\ZZ_p\backslash \{0\}$. Consider $\vx =
(w,\vx_p)\in\AA_N^\NN\times \PP_{\QQ_p}^1$ and define $\delta:=
d(({a\atop 0}),\vx_p)\cdot |a|_p$. Assume that $\delta\neq 0$. Let
$m\in\NN$ and $k\in\NN$ satisfy the inequality $p^k\ge
m\cdot(p\delta)^{-1}$. If~\eqref{eq_th_2} is satisfied:
$$
\{\vx_p, A^T \vx_p , \ldots, (A^T)^m \vx_p\} \subset B_k(\vx)
$$
then $\vx\not\in \lmad_\epsilon$ for $\epsilon= p\delta^{-1}m^{-1}$.

On the other hand if $\delta=0$ then $\vx\not\in\lmad$.
\end{theorem}

\section{Proof of Theorems \ref{th_main} and
\ref{th_da}}\label{sec_mainproof}

%

The distance $d(\cdot,\cdot)$ in $\PP_{\QQ_p}^1$ satisfies the
following property which will be widely used throughout the proof.
\begin{lemma}\label{fact_5}
Let $\ZZ_p$ be the set of $p$-adic integers and $A,B\in
\Sl^\pm_2(\ZZ_p)$ be two matrices such that $A\equiv B\pmod{p^k}$
for some $k\in\NN$, i.e. the $p$-adic norm of each entry of the
matrix $A-B$ is at most $p^{-k}$. Then for every two points
$\vw,\vv\in\PP_{\QQ_p}^1$ one has
$$
d(A\vw,\vv)\le \max\{d(B\vw,\vv), p^{-k}\}.
$$
\end{lemma}

\proof Denote $\vu=(u_1,u_2) = A\vw$ and $\vu' = (u_1',u_2') =
B\vw$. For simplicity we will choose $\vw,\vv\in
\tilde\PP_{\QQ_p}^1$. One can check that since $p$-adic norms of all
entries of $A$ are at most~1 then $\max\{|u_1|_p,|u_2|_p\}\le
\max\{|\omega_1|_p,|\omega_2|_p\}$. On the other hand since $A$ is
invertible then the inverse inequality is also true. This implies
that $\vu\in\tilde\PP^1_{\QQ_p}$ and by the same arguments
$\vu'\in\tilde\PP_{\QQ_p}^1$.

Now we calculate
$$
d(\vu,\vv) = |u_1\upsilon_2 - u_2\upsilon_1|_p = |u'_1\upsilon_2
-u'_2\upsilon_1 + d_1\upsilon_2 - d_2\upsilon_1|_p.
$$
where $\vd = (d_1,d_2) = (A-B)\vw$. Next,
$$
d(\vu,\vv)\le \max\{d(\vu',\vv), |d_1\upsilon_2 -
d_2\upsilon_1|_p\}.
$$
Since $A\equiv B\pmod{p^k}$ an upper bound for the second term of
the maximum is $p^{-k}$ which finishes the proof of the
lemma.\endproof

The next lemma shows that if some $\vq\in B_k(\vx)$ is not in
$\pbad_\delta$ then there exists some point $\vx_{p,l}$ which is not
in $\pbad_{\max\{\delta,p^{-k}\}}$.

\begin{lemma}\label{lem4}
Let $\vq=(q,q')\in B_k(\vx)$ and $\vu = (u,v)\in \ZZ^2$. Then if
$$|(\vu, \vq)|_p\cdot\min\{|q^{-1}|_p, |(q')^{-1}|_p\}\le \delta$$ for
some positive $\delta$ then there exists $l\in\NN$ such that
$$
|(\vu, \vx_{p,l})|_p\cdot
\min\{|x^{-1}_{p,l}|_p,|(x'_{p,l})^{-1}|_p\}\le \max\{\delta,
p^{-k}\}\quad\mbox{where}\quad \vx_{p,l} = (x_{p,l}, x'_{p,l}).
$$
\end{lemma}

{\bf Remark.} Here and afterwards by $(\vu,\vq)$ we denote the
standard inner product of vectors~$\vu$ and $\vq$.

\proof Since $\vq\in B_k(\vx)$ then there exists $\vx_{p,l}$ such
that $d(\vq, \vx_{p,l})\le p^{-k}$. Without loss of generality
assume that $|q'|_p\ge |q|_p$. Now we calculate
$$
|(\vu, \vx_{p,l})|_p = |ux_{p,l}+ vx'_{p,l}|_p =
|(q')^{-1}(uq'x_{p,l} - uqx'_{p,l} + uqx'_{p,l} +vq'x'_{p,l})|_p
$$$$
\le\max\{ \max\{|x_{p,l}|_p,|x'_{p,l}|_p\}\cdot |u|_p\cdot
d(\vq,\vx_{p,l}), |x'_{p,l}|_p\cdot |(q')^{-1}|_p\cdot|(\vu,
\vq)|_p\}
$$$$
\le \max\{|x_{p,l}|_p,|x'_{p,l}|_p\}\cdot \max\{p^{-k}, \delta\}.
$$
By dividing both sides of the inequality by
$\max\{|x_{p,l}|_p,|x'_{p,l}|_p\}$ we get the statement of the
lemma. \endproof

{\bf Proof of Theorem~\ref{th_main}.} The general idea of the proof
is to construct a sequence of points $\vq_n\in B_k(\vx)$, to show
that one of them do not lie in $\pbad_\epsilon$ and finally to apply
Lemma~\ref{lem4}. Define $\vq_0 =
(q_0,q_0')\in\tilde{\PP}_{\QQ_p}^1$ such that $\vq_0 = \vx_p$.
Represent the vector $\vq_0$ in the basis $\vw_1$ and $\vw_2$:
$\vq_0 = \alpha \vw_1 + \beta \vw_2$ where by Cramer's rule
$$
\alpha = \frac{\det (\vq_0,\vw_2)}{\det(\vw_1,\vw_2)},\quad \beta =
\frac{\det (\vw_1,\vq_0)}{\det(\vw_1,\vw_2)}.
$$
Since $\vq_0, \vx_1, \vw_2\in \tilde\PP_{\QQ_p}^1$ one can check
that
\begin{equation}\label{eq_ab}
|\alpha|_p =
\frac{d(\vq_0,\vw_2)}{d(\vw_1,\vw_2)}=\frac{\epsilon_2}{d(\vw_1,\vw_2)}
\quad\mbox{ and }\quad|\beta|_p =
\frac{d(\vq_0,\vw_1)}{d(\vw_1,\vw_2)}=\frac{\epsilon_1}{d(\vw_1,\vw_2)}.
\end{equation}
Define $\vq_n = (q_n,q'_n)$ as follows: $\vq_n = (A^T)^n \vq_0$.
Vectors $\vw_1$ and $\vw_2$ are eigenvectors of $A^T$ therefore
$$
\vq_n = \lambda_1^n\cdot \alpha \vw_1 + \lambda_2^n\cdot \beta
\vw_2.
$$
All entries of $(A^T)^n$ have $p$-adic norm at most 1, and moreover
$(A^T)^n$ is invertible. Therefore one can repeat the same arguments
as in the proof of Lemma~\ref{fact_5} to show that $\max\{|q_n|_p,
|q'_n|_p\}= \max\{|q_0|_p, |q'_0|_p\}=1$ or $\vq_n\in
\tilde\PP_{\QQ_p}^1$.

For an arbitrary integer vector $\vu = (u,v)$ one has $(\vu, \vq_n)
= A\lambda_1^n + B\lambda_2^n$ where $A = A(u,v) = \alpha\cdot(\vu,
\vw_1), B = B(u,v) = \beta\cdot(\vu,\vw_2). $

From the equation for eigenvalues $\lambda_1$ and $\lambda_2$ one
has $\lambda_1\cdot \lambda_2= \pm1$ and $|\lambda_1|_p =
|\lambda_{2}|_p = 1$. Therefore $\lambda_1^2 = \lambda_2^{-2}$. This
gives us the following
\begin{equation}
|\tilde{Q}_n|_p := |(\vu, \vq_{2\kappa(A)\cdot n})|_p = |A|_p \cdot
\left|\frac{B}{A} + \lambda^{n}\right|_p.
\end{equation}
where $\lambda = \lambda_1^{4\kappa(A)}$. Note that by
Lemma~\ref{fact_4}, $|\lambda-1|_p\le p^{-1}$ and therefore
$|\log\lambda|_p = |\log(\lambda_1^{4\kappa(A)})|_p = \epsilon_3\le
p^{-1}$.

In the rest of the proof we will construct integers $u$ and $v$ such
that for some $n\in \NN$ the value of $|\tilde{Q}_n|_p$ becomes so
small that $\vq_{2\kappa(A)\cdot n}$ is surely out of
$\pbad_\epsilon$. Then, using the fact that $\vq_{2\kappa(A)\cdot
n}$ is in $B_k(\vx)$ one concludes that $\vx$ is not in
$\lmad_\epsilon$.

We want to construct $u$ and $v$ such that $p\epsilon_3|A|_p>
|A+B|_p$. Then it will imply that $|1+B/A|_p<1$ and this will enable
us to use the idea from \cite{Bugeaud_Drmota_Mathan_2007}. By
Dirichlet pigeonhole principle one can choose
$(u,v)\in\ZZ^{2}\backslash \{0\}$ such that
$$
\begin{array}{l}
|A+B|_p = |uq_0+vq_0'|_p\le p(\epsilon_3\delta)^2;\\[4pt]
|u|,|v|< (\epsilon_3\delta)^{-1}.
\end{array}
$$
Note that since $v\in\ZZ$ then $|v|_p> \epsilon_3\delta$. Rewrite
the value $|(\vu,\vw_1)|_p$ in the following way
$$
|(\vu,\vw_1)|_p = |q_0|^{-1}_p\cdot | \omega_1(uq_o+vq_0') +
(\omega_2 q_0 - \omega_1 q_0')v|_p.
$$
%
Since $|\omega_1|_p\le 1$, an upper bound for the first summand is
$$
| \omega_1(uq_o+vq_0')|_p \le p(\epsilon_3\delta)^2
$$
and for the second summand it is
$$
|(\omega_2 q_0 - \omega_1 q_0')v|_p = \epsilon_1\cdot |v|_p >
\epsilon_3\delta^2\ge p(\epsilon_3\delta)^2.
$$
Therefore
$$
|(\vu,\vw_1)|_p = \frac{\epsilon_1\cdot |v|_p}{|q_0|_p}.
$$
By the construction we have $|A+B|_p \le p(\epsilon_3\delta)^2$ and
we just showed that
$$
|A|_p = |\alpha\cdot (\vu,\vw_1)|_p \stackrel{\eqref{eq_ab}}{=}
\frac{\epsilon_1\epsilon_2\cdot |v|_p}{|q_0|_p\cdot d(\vw_1,\vw_2)}.
$$
Notice that $|q_0|_p\cdot d(\vw_1,\vw_2) = |\omega_4
(q_0\omega_1-q_0'\omega_2)- \omega_2(q_0\omega_3-q_0'\omega_4)|_p\le
\max\{\epsilon_1,\epsilon_2\}$. This finally gives the following
lower estimate for $|A|_p$:
$$
\epsilon_3|A|_p\ge \epsilon_3\cdot
\min\{\epsilon_1,\epsilon_2\}\cdot |v|_p>(\epsilon_3\delta)^2\ge
p^{-1}\cdot |A+B|_p.
$$
So the aim is proved. Now we shall give an upper bound for $|A|_p$.
If $|q_0|_p<\epsilon\,$ where $\epsilon$ is given
by~\eqref{eq_epsilon} then we consider $\vu=(1,0)$ and
$|(\vu,\vq_0)|_p<\epsilon$ and by Lemma~\ref{lem4} with $\epsilon\ge
p^{-k}$ we get
$$
|(\vu, \vx_{p,l})|_p \cdot
\min\{|\vx_{p,l}^{-1}|_p,|(\vx'_{p,l})^{-1}|_p\}\le \epsilon
$$
for some $l\in\NN$. This implies that $\vx_{p,l}\not\in
\pbad_\epsilon$ and $(w,\vx_p)$ is not in $\lmad_{\epsilon}$ anyway.
Hence we can assume that $|q_0|_p\ge \epsilon$ and then
\begin{equation}\label{eq_ap}
|A|_p \le \epsilon_1\epsilon_2 \cdot \epsilon^{-1}\cdot
(d(\vw_1,\vw_2))^{-1}.
\end{equation}

Next, estimate the value $|\frac{B}{A}+\lambda^n|_p$. Since by
Lemma~\ref{fact_4}, $|\lambda - 1|_p\le p^{-1}$ and
$|\frac{B}{A}+1|_p<1$ we can use the property~\eqref{eq_log} of
$p$-adic logarithm:
$$
\left|\frac{B}{A} +\lambda^n\right|_p = \left|\log
\left(\frac{-B}{A}\right) - n\log \lambda\right|_p <
\left|\frac{\log(-B/A)}{\log\lambda} - n\right|_p.
$$
Now we show that $\log(-B/A)/\log\lambda$ lies in $\QQ_p$. This fact
is proved in~\cite{Bugeaud_Drmota_Mathan_2007} but for the sake of
completeness we repeat it here. It is trivial if $\lambda\in \QQ_p$.
Otherwise there exists a unique $\QQ_p$-automorphism $\sigma$ of
$\QQ_p(\lambda)$ different from the identity. We have
$\sigma(\lambda) = 1/\lambda$ and therefore $\sigma\log(\lambda) =
\log(\sigma(\lambda)) = -\log \lambda$. Next,
$$
A\lambda^n + B\lambda^{-n} = (\vu , \vq_{4\kappa(A)\cdot n}) =
\sigma ((\vu , \vq_{4\kappa(A)\cdot n})) = \sigma(A)\lambda^{-n} +
\sigma(B)\lambda^n.
$$
Since it is true for every natural $n$ then we have $\sigma(A) = B$
and $\sigma(B) = A$. Combining these equations for $\sigma$ together
we get
$$
\sigma\left(\frac{\log(-B/A)}{\log\lambda}\right) =
\frac{\log(-B/A)}{\log\lambda}.
$$
Further,
$$
\left|\log\left(-\frac{B}{A}\right)\right|_p =
\left|\frac{B}{A}+1\right|_p < p\epsilon_3 = p\cdot |\log\lambda|_p.
$$
therefore $|\log (-B/A)/\log\lambda|_p<p$. This together with the
fact that $\log (-B/A)/\log\lambda\in\QQ_p$ implies that this number
is also in $\ZZ_p$. Whence there exists integer $n$ within the range
$0\le n<p^d$ such that
$$
\left|\frac{\log(-B/A)}{\log \lambda} - n\right|_p\le p^{-d}.
$$
In other words there exists $n\in\ZZ$ within the range $0\le n<N$
such that
$$
\left|\frac{B}{A}+\lambda^n\right|\le \frac{p}{N}.
$$

One can choose $1\le n\le \frac{m}{2\kappa(A)}$ such that
$$
|\tilde{Q}_n|_p\le |A|_p\cdot \frac{2p\cdot
\kappa(A)}{m}\stackrel{\eqref{eq_ap}}\le
\frac{2\epsilon_1\epsilon_2\cdot p\cdot \kappa(A)}{\epsilon\cdot
d(\vw_1,\vw_2)\cdot m}.
$$
Then by substituting the formula~\eqref{eq_epsilon} for $\epsilon$
and using that $|u|,|v|< (\epsilon_3\delta)^{-1}$ we get
$\max\{u^2,v^2\}\cdot |\tilde{Q}_n|_p< \epsilon$. Since
$2\kappa(A)\cdot n\le m$, vector $\vq_{2\kappa(A)\cdot n}$ lies in
the set $B_k(\vx)$. Therefore Lemma~\ref{lem4} is applicable here
which in turn means that there exists $l\in \NN$ such that
$$
|(\vu, \vx_{p,l})|_p\cdot \min\{|x_{p,l}|_p^{-1},
|x'_{p,l}|_p^{-1}\}\le \max\{|\tilde{Q}_n|_p,p^{-k}\}.
$$
By~\eqref{eq_th_1} we estimate the value $p^{-k}$ to finally get the
bound
$$
|(\vu, \vx_{p,l})|_p\cdot \min\{|x_{p,l}|_p^{-1},
|x'_{p,l}|_p^{-1}\}<\epsilon\cdot \min\{u^{-2},v^{-2}\}.
$$
Hence $\vx_{p,l}\not\in \pbad_\epsilon$ and
$\vx\not\in\lmad_\epsilon$. \endproof

{\bf Proof of Theorem~\ref{th_da}.} As before represent $\vx_p$ as
$(q_0,q_0')$ where $\max\{|q_0|_p,|q_0'|_p\} = 1$ and denote it by
$\vq_0$. In this case one can easily check that $\delta =
|aq_0'|_p$. By $\vq_n$ we denote the point~$(A^T)^n \vq_0$ and as in
the proof of the previous theorem we have $\max\{|q_n|_p,|q_n'|_p\}
= 1$. Note that $(D_a^T)^n = D_{na}^T$ so one can easily get an
explicit formula for $\vq_n$:
$$
\vq_n = (q_0+na\cdot q_0', q_0').
$$

Firstly consider the situation when $\delta=0$. In that case we
should have $q_0'=0$ which straightforwardly implies that
$\vx_p\not\in \pbad_\epsilon$ for an arbitrarily small $\epsilon$.
One just need to consider $\vu = (1,0)$ and then $(\vu,\vq_0) = 0$.

Now we can assume that $\delta> 0$. By Dirichlet pigeonhole
principle one can choose $\vu = (u,v)\in\ZZ^2/\{0\}$ such that
$$
\begin{array}{l}
|(\vu,\vq_0)|_p\le p\delta^2;\\[4pt]
|u|,|v|< \delta^{-1}.
\end{array}
$$
Since $|u|_p> \delta$ and $|aq_0'|_p=\delta$ we get
$$
|(\vu,\vq_0)|_p\le p\delta^2 < p\cdot |u|_p\cdot |aq_0'|_p = p\cdot
|uaq_0'|_p.
$$
This implies that $|(\vu,\vq_0)|_p\le |uaq_0'|_p$ and therefore
$$
\frac{(\vu,\vq_0)}{uaq_0'}\in\ZZ_p
$$

One can write $(\vu,\vq_n)$ as
$$
(\vu,\vq_n) = (\vu,\vq_0) + nuaq_0',
$$
hence for each $d\in\NN$ there always exists an integer $n$ within
the range $0\le n<p^d$ such that
$$
|(\vu,\vq_n)|_p= |uaq_0'|_p\cdot \left|\frac{(\vu,\vq_0)}{uaq_0'} -
n\right|_p \le \delta p^{-d}.
$$
We choose the biggest possible $d$ such that $p^d\le m$ where $m$ is
defined in the statement of the theorem. Then we choose $0\le n<p^d$
such that $|(\vu,\vq_n)|_p\le \delta p^{-d}\le p\delta m^{-1}$ and
$$
\max\{u^2,v^2\}|(\vu,\vq_n)|_p < p\delta^{-1}m^{-1} = \epsilon.
$$
then since $n\le m$ we have $\vq_n\in B_k(\vx)$ and therefore by
Lemma~\ref{lem4} there exists $l\in\NN$ such that
$$
\max\{u^2,v^2\}\cdot |(\vu,\vx_{p,l})|_p\cdot
\min\{|x_{p,l}|^{-1}_p, |x'_{p,l}|^{-1}_p\}<\max\{u^2,v^2\}\cdot
\max\{|(\vu,\vq_n)|_p, p^{-k}\} \le \epsilon.
$$
Hence $\vx_{p,l}\not\in \pbad_\epsilon$ and $\vx\not\in
\lmad_\epsilon$.\endproof

\section{General restrictions on words imposed by Theorem~\ref{th_main}}

In this section we consider several corollaries from
Theorem~\ref{th_main} which will be more suitable for applications
to various particular words $w\in\AA_N^\NN$. They will also provide
some evidence that the conditions in this theorem together with
Theorem~\ref{th_da} are in fact quite restrictive.


Given $k\in \NN$ and an infinite or finite word $w$ over alphabet
$\AA_N$ we define the set
$$
U_k(w):=\{\phi_k(A_{w_n}))\in \Sl_2^\pm(\ZZ/p^k\ZZ)\; :\;n\in
\ZZ_{\ge 0}\}
$$
where $\phi_k\;:\Sl_2^\pm(\ZZ)\mapsto \Sl_2^\pm(\ZZ/p^k\ZZ)$ is the
canonical homomorphism. In other words it is the set of matrixes
correspondent to all the prefixes $w_n$ modulo $p^k$. Roughly
speaking Theorem~\ref{th_main} says that if $U_k(w)$ contains a
chain $A, A^2,\ldots, A^m$ with $A\in \widetilde{\Sl}_2(\ZZ_p)$ then
either $\vx_p$ coincides with one of the eigenvectors of $A$ or
$(w,\vx_p)$ is not in $\lmad_\epsilon$ for some $\epsilon$ which can
be explicitly calculated. Formally we have the following
\begin{corol}\label{corl2}
Let $\vx = (w,\vx_p)\in \AA_N^\NN \times \PP_{\QQ_p}^1$. Let
$\phi_k^*$ be a canonical homomorphism $\phi_k^*\;:\;
\Sl_2^\pm(\ZZ_p)\to \Sl_2^\pm(\ZZ/p^k\ZZ)$. Suppose that there
exists a matrix $A\in\widetilde{\Sl}_2(\ZZ_p)$ such that
\begin{equation}\label{eq_multA}
\phi^*_k(A), \phi^*_k(A^2),\ldots \phi^*_k(A^m)\in U_k(w)
\end{equation}
where the parameters $m,k\in\NN$ satisfy the
condition~\eqref{eq_th_1}. Assume also that the parameter $\delta$
defined in Theorem~\ref{th_main} for the matrix $A^T$ is not zero.
Then $\vx\not\in \lmad_\epsilon$ for $\epsilon$ given
by~\eqref{eq_epsilon}.
\end{corol}

To check it we basically use Lemma~\ref{fact_5} which shows
that~\eqref{eq_th_2} follows from~\eqref{eq_multA}. Then all the
conditions of Theorem~\ref{th_main} are satisfied.

In particular we can guarantee the condition~\eqref{eq_multA} if
$\phi^*_k(A)\cdot U_k(w) = U_k(w)$. Since $\Id\in U_k(w)$ then in
this case every integer power $\phi^*_k(A^m)$ will belong to
$U_k(w)$. In other words the condition~\eqref{eq_th_2} will be
always satisfied and the value $m$ from Corollary~\ref{corl2} will
be only restricted by~\eqref{eq_th_1}. As soon as $\vx_p$ does not
coincide with an eigenvector of $A$ one can make $m$ arbitrarily
large and $\epsilon$ arbitrarily small as $k$ tends to infinity. We
will show in a minute that the situation when $\phi^*_k(A)\cdot
U_k(w) = U_k(w)$ happens quite often.

One can easily check that $U_k(w)$ is always finite because the
space it lies in is finite. It is also easy to check that
$\#U_k(w)\to \infty$ as $k\to \infty$. The next proposition gives
information about the structure of sets $U_k(w)$.
\begin{proposition}\label{lem3}
Let $w$ be a recurrent word. Then for each $m\in \NN$ and $k\in \NN$
one has
\begin{equation}\label{eq_uk}
\phi_k(A_{w_m})\cdot U_k(\TT^m w) = U_k(w).
\end{equation}
\end{proposition}

\proof We first prove the following auxiliary statement: there
exists $m\in\NN$ such that $\phi_k(A_{w_m}) = \Id$. To show it we
construct the sequence $u_n$ of prefixes of $w$ by the following
rule.
\begin{itemize}
\item $u_1:= w_1$.
\item Given $u_n$ we take the prefix $u_{n+1}$ such that it ends
with $u_n$. In other words $u_{n+1} = v_{n+1}u_n$ for some word
$v_{n+1}$. We can always do that because $w$ is recurrent.
\end{itemize}
Because the set $\Sl_2^\pm(\ZZ/p^k\ZZ)$ is finite we can find
$s,t\in\NN$ such that $\phi_k(A_{u_s})=\phi_k(A_{u_{s+t}})$. Then
notice that
$$
u_{s+t} = v_{s+t}v_{s+t-1}\ldots v_{s+1}u_s.
$$
By substituting this to the matrix equality we get that
$\phi_k(A_{v_{s+t}\ldots v_{s+1}}) = \Id$. Finally the observation
that $v_{s+t}\ldots v_{s+1}=w_m$ is the prefix of $w$ finishes the
proof of the auxiliary statement.

Now we can prove the proposition. We will prove~\eqref{eq_uk} for
$m=1$. The rest can easily be done by induction. The inclusion
$\phi_k(A_{w_1})\cdot U_k(\TT w) \subseteq U_k(w)$ is
straightforward. Indeed for every prefix $u$ of $\TT w$ the word
$w_1u$ is the prefix of $w$ and therefore $\phi_k(A_{w_1}A_u) =
\phi_k(A_{w_1u})\in U_k(w)$. For inverse inclusion we only need to
check that $\Id$ belongs to $\phi_k(A_{w_1})\cdot U_k(\TT w)$. Other
elements of $U_k(w)$ correspond to prefixes $w_n$ which start with
the letter $w_1$ and this fact straightforwardly implies that they
also belong to $\phi_k(A_{w_1})\cdot U_k(\TT w)$. On the other hand
by the auxiliary statement one can represent $\Id$ as
$\phi_k(A_{w_m})$ where the word $w_m$ starts with~$w_1$ as well.
The proof is finished.
\endproof

Since there are only $\#U_k(w)$ different elements $\phi_k(A_{w_m})$
then Proposition~\ref{lem3} shows that the collection $\UU_k:= \{
U_k(\TT^m w)\;:\; m\in \ZZ_{\ge 0}\}$ consists of at most $\#U_k(w)$
elements. On the other hand if $\#\UU_k <\#U_k(w)$ then one can find
two prefixes $w_m, w_n$ with $\phi_k(A_{w_m})\neq \phi_k(A_{w_n})$
such that
$$
\phi_k(A^{-1}_{w_m})\cdot U_k(w) =\phi_k(A^{-1}_{w_n})\cdot U_k(w).
$$
Without loss of generality assume that $m>n$. Then this leads us to
$\phi_k(A_u)\cdot U_k(\TT^n w) = U_k(\TT^n w)$ where $w_nu = w_m$.
By Lemma~\ref{fact_6}, $A_u\in \widetilde{\Sl}_2(\ZZ_p)$ and we also
have $\phi_k(A_u)\neq \Id$ so the condition~\eqref{eq_multA} is
satisfied for $A = A_u$. Then the application of
Corollary~\ref{corl2} will give us that either $\vx_p$ is an
eigenvalue of $A_u$ or $\TT^n\vx\not \in \lmad_\epsilon$ for
$\epsilon$ defined by~\eqref{eq_epsilon}. However with help of only
this basic observation we do not have a big control on $\epsilon$.
To use Corollary~\ref{corl2} in full we need some quantitative
result of this kind which will be true for an arbitrarily large~$k$.

Another useful corollary of Proposition~\ref{lem3} is that for every
$m\in\NN$ the sets $U_k(\TT^mw)$ have the same size.

 For every $n\in\NN$ and every infinite word $w\in\AA_N^\NN$ we define the
 sets
$$
V(n,w):=\{A_{w_m}\;:\; w_n\mbox{ is the prefix of }\TT^mw \}
$$
and $V_k(n,w):= \phi_k(V(n,w))$. One can easily check that the sets
$V(n,w)$ are nested as~$n$ grows: $V(n_1,w)\subset V(n_2,w)$ for
$n_1\ge n_2$. The same is surely true for their projections
$V_k(n,w)$. If $w$ is recurrent then $V(n,w)$ is always infinite for
all $n\in\NN$, however $V_k(n,w)$ is obviously finite. We also
define the set $Vp(n,w)$ which coincides with $V(n,w)$ but every
element in $Vp(n,w)$ is considered as a matrix over $\ZZ_p$. Finally
let
$$
Vp(w):= \bigcap_{n=1}^\infty \overline{Vp(n,w)}.
$$
We can not say anything about whether it is finite or not. However
it is always non-empty since obviously $\Id\in Vp(w)$.

\begin{proposition}\label{prop1}
Let $w\in \AA_N^{\NN}$ be recurrent word. For each $k\in \NN$ there
exists $n\in\NN$ such that for every $A\in V_k(n,w)$ we have $A\cdot
U_k(w) = U_k(w)$.
\end{proposition}


\proof For each matrix $B\in U_k(w)$ associate the value $n(B)$, the
length of the shortest word $w_{n(B)}$ such that $B =
\phi_k(A_{w_{n(B)}})$. Since there are finitely many elements in
$U_k(w)$ we can define $n:= \max_{B\in U_k(w)} \{n(B)\}$. By the
construction of $n$ for every infinite word $w^*$ which starts with
$w_n$ one has $U_k(w^*)\supset U_k(w)$.

Now for each $m\in\NN$ such that $A_{w_m}\in V(n,w)$ by definition
we have that the word $\TT^m w$ starts with $w_n$. Therefore an
application of Proposition~\ref{lem3} gives us $U_k(w)\subset
U_k(\TT^m w) = \phi_k(A_{w_m})\cdot U_k(w)$. Or in other words,
$\forall A\in V_k(n,w)$ one has $A\cdot U_k(w) \supset U_k(w)$.
Finally since two sets $U_k(w)$ and $U_k(\TT^mw)$ have the same
cardinality we have the equality $A\cdot U_k(w) = U_k(w)$.
\endproof

\section{Relation between $Vp(w)$ and $\lmad$. Proof of
Theorem~\ref{th_concat}}\label{sec_concat}

Proposition~\ref{prop1} shows that the condition~\eqref{eq_multA} is
satisfied for each $A\in V_k(n,w)$ where $n$ is large enough. This
observation gives rise to the following

\begin{theorem}\label{th_plc_cond}
Let $A_1,A_2\in \widetilde{\Sl}_2(\ZZ_p)$ be such that their
eigenvectors $\vw_1,\vw_2,\vw_3,\vw_4\in\PP_{\overline{\QQ}_p}^1$
are all distinct. Let $w\in \AA^\NN_N$. Assume that for some
$k\in\NN$ and for all $n\in \NN$, $\phi_k(A_1), \phi_k(A_2)\in
V_k(n,w)$. Then for each $\vx_p\in \PP_{\QQ_p}^1$, $(w,\vx_p)\not\in
\lmad_\epsilon$ where $\epsilon \gg_{A_1,A_2,p} p^{-k}$.
\end{theorem}

Here we used the Vinogradov symbol: $x\gg_{A_1,A_2,p} y$ means that
$x\ge c\cdot y$ where $c$ is a positive constant dependent on
$A_1,A_2$ and $p$ only.

\proof Consider the point $\vx = (w,\vx_p)$. Since all the values
$\vw_1,\ldots,\vw_4$ are distinct then
$$
\max\{ \min\{d(\vx_p,\vw_1),d(\vx_p,\vw_2)\},
\min\{d(\vx_p,\vw_3),d(\vx_p,\vw_4)\}\} \gg_{A_1,A_2,p} 1.
$$
Without loss of generality assume that the maximum above is reached
for $\vw_1$ and $\vw_2$. Choose the matrix $A = A_1$. Then the
estimate on the parameter $\delta$ from Theorem~\ref{th_main} is
$\delta\gg_{A_1,A_2,p} 1$.

We have that $\phi_k(A)$ is in $V_k(n,w)$ for all $n\in\NN$.
Therefore by Proposition~\ref{prop1} we have $\phi_k(A)\cdot U_k(w)
= U_k(w)$ and therefore \eqref{eq_multA} is true for an arbitrary
$m\in\NN$. Choose the biggest possible $m$ such that the
condition~\eqref{eq_th_1} is satisfied, then $m\asymp p^{2k}$. Now
Corollary~\ref{corl2} implies that
$$
\vx\not\in\lmad_\epsilon\;\mbox{ for }\; \epsilon
\stackrel{\eqref{eq_epsilon}}\gg_{A_1,A_2,p} m^{-1/2}\gg_{A_1,A_2,p}
p^{-k}.
$$
\endproof

Note that in Theorem~\ref{th_plc_cond} the estimate for $\epsilon$
depends only on $A_1,A_2$ and $p$ but it does not depend on $k$.
This simple observation makes the following corollary true.

\begin{corol}\label{corl_plc_cond}
If in terms of Theorem~\ref{th_plc_cond} one additionally has that
$A_1,A_2\in Vp(w)$ then $\vx\not\in\lmad$.
\end{corol}

Corollary shows that for $(w,\vx_p)\in\lmad$ the set $Vp(w)$ must be
very small in a sense that it should not contain two ``independent''
matrices, i.e. the matrices from $\widetilde{\Sl}_2(\ZZ_p)$ which
eigenvectors do not intersect. This observation is sufficient to
prove Theorem~\ref{th_concat}.

{\bf Proof of Theorem~\ref{th_concat}.} Firstly note that for $n> m$
the word $\sigma_{n-1}$ is the prefix of $\sigma_n$ and therefore it
is the prefix of $\sigma_l$ for any $l\ge n$. This in turn implies
that $\sigma_{n-1}$ is a prefix of $w$. Let $L(n)$ be the length
of~$\sigma_n$. Then for $n>2m$, $\sigma_{n-m}$ is the prefix of
$\TT^{L(n-1)}\sigma_n$. These two facts together imply that
$A_{\sigma_{n-1}}\in V(L(n-m),w)$.

Now fix $k\in\NN$ and consider the sequence $s_n :=
\{\phi_k(A_{\sigma_n})\}_{n\in\NN}$. Since $s_{n+m}$ depends only on
$m$-tuple $s_n,\ldots,s_{n+m-1}$ and there are only finitely many of
such $m$-tuples over $\Sl_2^\pm (\ZZ/p^k\ZZ)$ then the sequence
$s_n$ is eventually periodic. Moreover by the condition of the
theorem $s_n$ is also uniquely determined by
$s_{n+1},\ldots,s_{n+m}$ therefore $s_n$ is purely periodic.

$L(n)$ tends to infinity as $n\to\infty$ therefore $s_1,\ldots,s_m$
are in $V_k(l,w)$ for an arbitrary large~$l$. Then since the
inclusion is true for all $k\in\NN$ we have that $s_1,\ldots,s_m\in
Vp(l,w)$ for all $l\in\NN$ and finally
$A_{\sigma_1},\ldots,A_{\sigma_k}\in Vp(w)$. Choose two of these
matrices which are different (by the condition of the theorem we can
do so). Without loss of generality let them be $A_{\sigma_1}$ and
$A_{\sigma_2}$. Since $\sigma_1$ and $\sigma_2$ are one letter words
we use Lemma~\ref{fact_6} to show that $A_{\sigma_1}$ and
$A_{\sigma_2}$ are both in  $\widetilde{\Sl}_2(\QQ_p)$ and that
their eigenvectors $\vw_{i}$, $i\in\{1,2,3,4\}$ are all distinct.
Therefore Corollary~\ref{corl_plc_cond} can be applied to $\vx =
(w,\vx_p)$ which finishes the proof of the theorem.\endproof

We finish this section by showing that even weaker condition than in
Corollary~\ref{corl_plc_cond} on $Vp(w)$ tells quite a lot about
possible elements $(w,\vx_p)\in\lmad$.

\begin{theorem}\label{th_vpw}
Assume that $A\in\widetilde{\Sl}_2(\ZZ_p)$ is an element of $Vp(w)$
and $\vw_1,\vw_2$ are two eigenvectors of $A$. Then there are at
most two different values $\vx_p$ such that $\vx=(w,\vx_p)\in\lmad:$
$\vx_p = \vw_1$ or $\vx_p = \vw_2$.
\end{theorem}

Note that both $\vw_1$ and $\vw_2$ are quadratic algebraic numbers
therefore if $w$ satisfies the conditions of Theorem~\ref{th_vpw}
and $(w,\vx_p)\in\lmad$ then $\vx_p$ must be quadratic irrational.

\proof there are two possible values $\vx_p$ for which $\delta$
defined in Theorem~\ref{th_main} is zero: $\vx_p = \vw_1$ and $\vx_p
= \vw_2$ where $\vw_1$ and $\vw_2$ are eigenvectors of $A$. If
$\vx_p$ is not one of that two values then~$\delta$ becomes strictly
positive and we can repeat the proof of Theorem~\ref{th_plc_cond}
with $A$ in place of $A_1$ to show that $(w,\vx_p)\not \in \lmad$.
\endproof

In view of Corollary~\ref{corl_plc_cond} and Theorem~\ref{th_vpw} it
would be interesting to investigate sets $Vp(w)$ for various words
$w\in\AA_n^\NN$. In particular it would be good to describe the
collection of words for which $Vp(w)$ does not contain any matrices
$A\in\widetilde{\Sl}_2(\ZZ_p)$.

\section{Periodic words $w$}\label{sec_5}


Now we prove Theorem~\ref{th_lmad}. Consider an arbitrary purely
periodic word $w$. Let $u$ be the finite factor of $w$ which
comprises the minimal period of $w$. Then one can check that $Vp(w)
= \{A_u^k\}_{k\in\NN}$. By periodicity $\TT^l w = w$ where $l$ is
the length of $u$.

We start with more difficult ``if'' part of the theorem. In that
case we are given that $\vx_p$ coincides with an eigenvector of
$A_u^T$. In other words it is the root of equation
\begin{equation}\label{eq_quad_q}
\TT^l(w,\vx_p) = (w,\vx_p)
\end{equation}
If we represent $\vx_p$ in affine form $\vx_p = ({x_p\atop 1})$ then
it becomes a quadratic equation $a+cx_p = (b+dx_p)x_p$ where $A_u =
({a\; b\atop c\; d})$.

\begin{lemma}\label{lem_quad_ir}
Let $v\in\QQ_p\backslash \QQ$ be the solution of quadratic equation
with integer coefficients. Then $({v\atop 1})\in \pbad_\epsilon$ for
some $\epsilon = \epsilon(v)>0$.
\end{lemma}

The proof of the lemma uses similar ideas as Liouville's proof that
every real quadratic irrational is badly approximable.

\proof Let $\bar v$ be a conjugate of $v$. Consider $a,b\in\ZZ$ such
that $|av-b|_p<1$, for other pairs $a,b$ the conditions of
$\pbad_\epsilon$ are held automatically. Since $|a(v-\bar v)|_p\ll_v
1$ we have that $|a\bar v - b|_p \ll_v 1$.

Consider the value $ (av-b)(a\bar v-b)$. It is integer, therefore
$$
|(av-b)(a\bar v-b)|_p\ge (a^2 v\bar v-ab (v+\bar v)+b^2)^{-1} \gg_v
(\max\{a^2,b^2\})^{-1}.
$$
Hence we get
$$
|av-b|_p\cdot \min\{1,|v^{-1}|_p\} = \frac{|(av-b)(a\bar
v-b)|_p\cdot \min\{1,|v^{-1}|_p\}}{|a\bar v-b|_p}\gg_v
(\max\{a^2,b^2\})^{-1}.
$$ \endproof

\begin{lemma}\label{lem_conform}
Let $\vx\in\pbad$. Then for any matrix $A\in \Gl_2(\ZZ)$ the point
$A \vx$ is in $\pbad$.
\end{lemma}

\proof Choose $\vx\in \tilde\PP_{\QQ_p}^1$, i.e. $\vx=(x,y)$ with
$\max\{|x|_p,|y|_p\}=1$. Then since $\vx\in\pbad$ we have
$|ux+vy|_p\gg \min\{|u|^{-2},|v|^{-2}\}$.

Write $A$ in coordinate form: $A = ({a\;b\atop c\;d})$. Consider
$|u(ax+by)+v(cx+dy)|_p$. Then
$$
|u\cdot (ax+by)+ v\cdot(cx+dy)|_p =|(ua+vc)x + (ub + vd)y|_p $$$$\gg
\min\{|ua+vc|^{-2}, |ub+vd|^{-2}\}\gg \min\{|u|^{-2},|v|^{-2}\}.
$$ \endproof

Now we are ready to show that for any solution of~\eqref{eq_quad_q}
we have $(w,\vx_p)\in\lmad$. From the theory of continued fractions
we know that $x_p\not\in \QQ$ (the solution of the equation $a+cx_p
= (b+dx_p)x_p$ over $\RR$ gives a solution which continued fraction
expansion is periodic with period~$u$, so it is infinite). Therefore
by Lemma~\ref{lem_quad_ir}, $\vx_p$ is in $ \pbad_{\epsilon_0}$ for
some positive~$\epsilon_0$. Then by the construction of $w$ and
$\vx_p$ we have that $\TT^l (w,\vx_p) = (w,\vx_p)$ therefore to
ensure that $(w,\vx_p)\in \lmad$ it is sufficient to check that
$\vx_{p,i}\in \pbad_{\epsilon_i}$ where $1\le i< l$, $ \TT^i
(w,\vx_p) = (\TT^i w, \vx_{p,i})$ and $\epsilon_i$ are some positive
constants. However this is in fact true by Lemma~\ref{lem_conform}
and the fact that by the construction of~$\TT$, $\vx_{p,i} =
(A_{w_i})^T\vx_p$. So, $(w,\vx_p)\in\lmad_{\epsilon}$ where
$\epsilon = \min_{0\le i< l}\{\epsilon_i\}$. Now note that if
$\TT^l(w,\vx_p) = (w,\vx_p)$ then $\vx_p$ coincides with one of the
eigenvectors of $A^T_{w_l}$. This finishes ``if'' part of the proof
of Theorem~\ref{th_lmad}.

On the other hand the ``only if'' part this theorem is a
straightforward corollary of Theorem~\ref{th_main}. We just take $A
= A_{w_l}$ which by Lemma~\ref{fact_6} is in
$\widetilde{\Sl}_2(\ZZ_p)$. If $\vx_p$ does not coincide with any of
eigenvectors of $A^T$ then the parameter $\delta$
in~Theorem~\ref{th_main} is non-zero. Finally take $k$ arbitrarily
large and take the biggest $m\in\NN$ which
satisfies~\eqref{eq_th_1}.

\section{Collections $\UU_k$ and low complexity words
$w$}\label{sec_lowcomplex}

In this section we show how can Theorems~\ref{th_main}
and~\ref{th_da} can be applied to words with low complexity. We
start with the following auxiliary lemma.

\begin{lemma}\label{lem9}
Let $u,v$ be two finite words over the alphabet $\NN$. If $A_u,
A_v\in\widetilde{\Sl}_2(\ZZ_p)$ share the same eigenvector $\vv$
then there exists a finite word $w$ and positive integer values
$m_1,m_2$ such that $u = w^{m_1}, v = w^{m_2}$. Moreover $A_w$
shares the same eigenvector $\vv$.
\end{lemma}

\proof We consider $\vv$ in affine form, so we can look at it as a
number from $\overline{\QQ}_p$. Note that $\vv$ is algebraic and,
since there is a field isomorphism between algebraic numbers in
$\overline{\QQ}_p$ and in $\CC$, two matrices $A_u$ and $A_v$ should
also share the same eigenvector $\vv^*\in \PP_\CC^1$.

From the theory of continued fractions we know that $A_u$ always has
two different real eigenvectors. Considered as real numbers they
satisfy the conditions $\vv^*_1<0<\vv_2^*$ and $w_{CF}(\vv_2^*) =
u^\infty \,(=uuuuuu....)$. Since $A_u$ and $A_v$ share the same
eigenvector in $\PP_{\QQ_p}^1$ they must also share the same
positive eigenvector in $\PP_\CC^1$ which in turn implies that two
infinite words $u^\infty$ and $v^\infty$ coincide. The conclusion of
the lemma can be easily derived from this fact.\endproof

Now we are ready to prove
\begin{proposition}\label{lem6}
Let $w\in\AA_N^\NN$ be recurrent and $\UU_k$ be constructed from
$w$. If for every $k\in\NN$, $\#\UU_k$ is bounded above by some
absolute constant then there are at most two points $\vx_p\in
\PP_{\QQ_p}^1$ such that $\vx = (w,\vx_p)\in \lmad$. Moreover if $w$
is not periodic then $(w,\vx_p)\not\in \lmad$ for all points
$\vx_p\in \PP_{\QQ_p}^1$.
\end{proposition}

\proof Assume that the sequence $\#\UU_k$ is bounded. Since
$\#\UU_k$ is non decreasing then there exists $k_0$ such that for
$k\ge k_0$, $\#\UU_k$ is a constant. Consider the minimal number
$n\in\NN$ such that $\phi_{k_0}(A_{w_n}) = \Id$, from the proof of
Proposition~\ref{lem3} we know that such $n$ always exists. Then by
Proposition~\ref{lem3}, $\forall k\in\NN$, $U_k(w) =
\phi_k(A_{w_n})\cdot U_k(\TT^n w)$ and in particular $U_{k_0}(w) =
U_{k_0}(\TT^n w)$. Remind that
$$
\UU_k := \{U_k(\TT^n w)\;:\; n\in\NN\}.
$$
Since for every $k\ge k_0$ the number of elements in $\UU_k$ stays
the same we should have $U_k(w) = U_k (\TT^n w)$ for every
$k\in\NN$. In particular it means that $\forall k\in\NN$,
$$
\phi_k(A_{w_n})\cdot U_k(w) = U_k(w)
$$
and the conditions of Corollary~\ref{corl2} are satisfied for $A =
A_{w_n}$. As $k$ tends to infinity one can choose an arbitrary large
$m$ satisfying~\eqref{eq_th_1}. If $\vx_p$ does not coincide with
any of the eigenvectors of $A^T$ then we have $\delta>0$, so
Corollary~\ref{corl2} can be applied to get $\vx\not \in \lmad$.
This shows the first statement of the proposition.

For the second statement we consider an infinite sequence
$n_1<n_2<\ldots<n_t<\ldots$ of positive integers such that
$\phi_{k_0}(A_{w_{n_t}}) = \Id$. By slightly modifying the arguments
of Proposition~\ref{lem3} one can show that such a sequence also
exists. We showed that if $(w,\vx_p)\in \lmad$ then $\vx_p$ must
coincide with one of the eigenvalues of all matrices $A_{w_{n_t}}$,
in other words all of them must share the same eigenvalue. By
Lemma~\ref{lem9} it means that there exists a finite word
$\tilde{w}$ such that $w_{n_t} = \tilde{w}^{m_t}$ which
straightforwardly implies that $w$ is periodic.
\endproof

We will associate each word $w\in\AA_N^\NN$ with another word $u =
u(k)\in\UU_k^\NN$ by the following rule: $u = b_1b_2\ldots$ where
$$
b_n = U_k(\TT^{n-1} w).
$$

\begin{lemma}\label{lem7}
Let $\vx=(w,\vx_p)\in \lmad$. Then there exists $k_0\in\NN$ such
that for every $k>k_0$ and every two-letter factor $c_1c_2$ of
$u(k)$ there is the unique value $a\in \AA_N$ with the following
property:
\begin{itemize}
\item if $b_nb_{n+1} = c_1c_2$ then $a_n = a$ where $a_n$ is the
$n$th letter of $w$.
\end{itemize}
\end{lemma}

In other words Lemma~\ref{lem7} states that for $k$ large enough any
two-letter factor of $u$ uniquely determines a one-letter factor of
$w$.

\proof Assume the contrary: one can find an arbitrarily large
$k\in\NN$ such that there exist positive integer $n$ and $l$ such
that $b_nb_{n+1} = b_lb_{l+1}=c_1c_2$ but the corresponding letters
$a_n = a$ and $a_l = a'$ are different. From Proposition~\ref{lem3}
we have
$$
c_1 = \phi_k(A_a)\cdot c_2;\quad\mbox{and}\quad c_1 =
\phi_k(A_{a'})\cdot c_2.
$$
which immediately implies that
$$c_1 = \phi_k(A_{a'}^{-1}A_a)\cdot c_1 = \phi_k(D_{a-a'})c_1.$$
In other words it means that $U_k(\TT^n w)$ contains matrices
$$
\phi_k(D_{a'-a}), \phi_k(D_{a'-a}^2), \ldots,
\phi_k(D_{a'-a}^m),\ldots
$$
Lemma~\ref{fact_5} shows that this property implies that
$$
\{\vx_{p,n}, D_{a-a'}^T \vx_{p,n},\,\ldots\,, (D_{a-a'}^T)^m
\vx_{p,n}\}\subset B_k(\TT^n\vx)
$$
where $m$ can be made arbitrarily large. In other
words~\eqref{eq_th_2} is satisfied for $\TT^n\vx$. Take $m =
p^{k+1}\delta$. Then all the conditions of Theorem~\ref{th_da} are
satisfied for the point $\TT^n\vx$. If $\vx_{p,n}$ is such that
$\delta=0$ then Theorem~\ref{th_da} states that
$\TT^n\vx\not\in\lmad$ which by the invariance of $\lmad$ implies
that $\vx\not\in\lmad$. This contradicts the conditions of the
lemma. Otherwise it implies that $\TT^n \vx \not \in \lmad_\epsilon$
where $\epsilon$ can be made an arbitrarily small positive number as
$k$ grows to infinity. Again this leads to a contradiction.
\endproof

In particular Lemma~\ref{lem7} shows that if $\vx\in\lmad$ and the
word $u$ is periodic then the word $w$ must be periodic too.
Theorem~\ref{th_lmad} gives a complete description of elements
$(w,\vx_p)\in\lmad$ in this case.

Now we will check that for $l$ large enough $l$-letter factor of $w$
uniquely determines $(l+1)$-letter factor of $u$.

\begin{lemma}\label{lem8}
Let $w$ be uniformly recurrent word. There exists $l\in\NN$ such
that for every $l$-letter factor $v$ of $w$ there exists unique
$l+1$-letter factor $v'$ of $u$ such that
\begin{itemize}
\item if $v = (\TT^n w)_l$ then $v' = (\TT^n u)_{l+1}$.
\end{itemize}
\end{lemma}

\proof From the proof of Proposition~\ref{prop1} we know that there
exists $m\in\NN$ such that for every infinite word $w^*$ starting
with $w_m$ we have $U_k(w)\subset U_k(w^*)$. Since $w$ is uniformly
recurrent there exists an absolute constant $m_1$ such that the
distance between two consecutive factors $w_m$ in $w$ is at most
$m_1$. Then one can easily check that every factor of $w$ of length
$l = m+m_1$ contains $w_m$ as a factor. We will show that this value
$l$ works for the lemma.

Consider an arbitrary factor $v$ of $w$ of length $l$. We showed
that it can be written as $v = u_1w_mu_2$. Denote the length of
$u_1$ by $n_1$. Consider an arbitrary $n\in\NN$ such that $(\TT^n
w)_l = v$. Then $(\TT^{n+n_1}w)_m = w_m$ and $U_k(\TT^{n+n_1}
w)\supset U_k(w)$. However every set in the collection~$\UU_k$ has
the same cardinality therefore $U_k(\TT^{n+n_1} w) = U_k(w)$.
Finally by Proposition~\ref{lem3}, $U_k(\TT^n w) =
\phi_k(A_{v_1})\cdot U_k(w)$. The right hand side of the last
equality does not depend on the position~$n$ of the factor $v$ so
the lemma is proved.
\endproof


{\bf Proof of Theorem~\ref{th_sturmian}.} Firstly without loss of
generality we can assume that $w$ is uniformly recurrent. Otherwise
$\overline{\{\TT^n \vx\}}_{n\in\NN}$ is not minimal set invariant
under $\TT$. Then we can consider its minimal subset which will be
of the form $\overline{\{\TT^n(w',\vx_p')\}}_{n\in\NN}$ where $w'$
is uniformly recurrent. Since $P(w',n)\le P(w,n)$ all the conditions
of Theorem~\ref{th_sturmian} will be satisfied for $w'$ as well and
the statement of the theorem for $w$ will easily be followed from
the same statement for $w'$.

Suppose that a value $\vx_p\in\PP_{\QQ_p}^1$ such that
$(w,\vx_p)\in\lmad$ exists. If for some $n\in\NN$ we have that
$P(w,n+1) = P(w,n)$ then $w$ is periodic (see~\cite{Lothaire_2001}
for details). But this contradicts to the conditions of the theorem.
So we should have $P(w,n+1) - P(w,n)\ge 1$ for every $n\in\NN$. Also
Proposition~\ref{lem6} shows that the sequence
$\{\#\UU_k\}_{k\in\NN}$ is unbounded.

Since $P(w,n)\le n+C$, there exists a value $n_0\in \NN$ such that
for every $n\ge n_0$, $P(w,n+1) - P(w,n) = 1$. Choose $k_0$ such
that $\#\UU_{k_0} > P(w,n_0)$. This will ensure that for each
$k>k_0$, $P(w,n_0)<P(u(k),n_0+1)$. If for some $n\in \NN$, $P(u,n) =
P(u,n+1)$ then $u$ is periodic which by Lemma~\ref{lem7} implies
that $w$ is periodic too --- a contradiction. Hence we have
$P(w,n+1) - P(u,n)\ge 1$ and this immediately implies that $\forall
n\ge n_0$, $P(w,n)<P(u,n+1)$. On the other hand Lemmata~\ref{lem7}
and~\ref{lem8} together imply that there exists $l_0\in\NN$ such
that $\forall l\ge l_0$, $P(u,l) = p(w,l-1)$. So we get a
contradiction with assumption that the value $\vx_p$ exists.
\endproof

\noindent{\bf Remark.} The complexity condition in
Theorem~\ref{th_sturmian} does not seem to be sharp. One can
possibly use more delicate arguments to improve them.

\section{Proof of Theorem~\ref{th_graph}}\label{sec_graph}

%


Since $w$ is uniformly recurrent, Lemma~\ref{lem8} can be applied to
find the value $l=l(k)\in\NN$ such that every factor $v$ of $w$ of
length $l$ uniquely determines an element $u\in\UU_k$. In other
words as soon as $\TT^nw$ starts with $v$ the set $U_k(\TT^nw)$
remains the same. Moreover the proof of Lemma~\ref{lem8} tells that
$U_k(\TT^nw) = U_k(v)$.

Now consider the graph $G_l$ and take two vertices $t_1,t_2\in T_l$
such that they are connected with the same vertex $s\in S_l$. then
$U_k(s) = U_k(st_1) = U_k(st_2)$. We apply Proposition~\ref{lem3} to
get
$$U_k(s) = \phi_k(A_s) U_k(t_1) = \phi_k(A_s)
U_k(t_2)\quad\mbox{or}\quad U_k(t_1) = U_k(t_2).
$$
By repeating this argument for every triple $t_1,t_2,s$ in the graph
$G_l$ one can show that for any two words $t^*_1, t^*_2\in T_l$ from
the same connected component of $G_l$ we must have $U_k(t^*_1) =
U_k(t^*_2)$. The number of different elements $U_k(s)$ where $s$
runs through all vertices in $S_l$ coincides with~$\#\UU_k$.
Therefore if the number of connected components in $G_l$ is bounded
by absolute constant independent on $l$ then $\#\UU_k$ is bounded by
the same constant. However Proposition~\ref{lem6} states that in
this case $(w,\vx_p)\not\in\lmad$ for every $\vx_p\in
\PP_{\QQ_p}^1$.
\endproof

It would be interesting to investigate the graphs $G_l$ for various
infinite words $w$. They are surely connected for infinite words $w$
such that $\{\TT^nw\}_{n\in\NN}$ is dense everywhere on $\AA^\NN_N$.
On the other hand numerical evidence suggests that the number of
connected components of $G_l$ for Thue-Morse word $w_{tm}$ tends to
infinity as $l\to\infty$. Even though the author does not know the
proof of this fact it seems that Thue-Morse word is not covered by
the last theorem.

\vspace{5mm}

\noindent Dzmitry A. Badziahin: Department of Mathematics, Durham
University,

\vspace{0mm}

\noindent\phantom{Dzmitry A. Badziahin: }Durham, DH1 3LE, England.


\noindent\phantom{Dzmitry A. Badziahin: }e-mail:
dzmitry.badziahin@durham.ac.uk


\begin{thebibliography}{99}

\bibitem{BBEK_2014}
D. Badziahin, Y. Bugeaud, M. Einsiedler, D. Kleinbock, On the
complexity of a putative counterexample to the $p$-adic Littlewood
conjecture. Preprint, 2014.

\bibitem{Bugeaud_Drmota_Mathan_2007}
Y. Bugeaud, M. Drmota, and B. de Mathan, On a mixed Littlewood
conjecture in Diophantine approximation. {\it Acta Arith.} V. 128
pp. 107--124, 2007.

\bibitem{Cassels_1955}
J. W. S. Cassels, An introduction to Diophantine Approximation.
Cambridge University Press, Cambridge, 1955.

\bibitem{Einsiedler_Kleinbock_2007}
M. Einsiedler and D. Kleinbock, Measure rigidity and $p$-adic
Littlewood-type problems. {\it Compositio Math.} V. 143 pp.
689--702, 2007.

\bibitem{Lothaire_2001}
M. Lothaire, Algebraic combinatorics on words. Cambridge University
Press, Cambridge, 2002.

\bibitem{Mathan_Teulie_2004}
B. de Mathan et O. Teuli\'e, Probl\`emes diophantiens simultan\'es.
{\it Monatsh. Math.} V. 143 pp. 229--245, 2004.


\end{thebibliography}
\end{document}